\documentclass[11pt,a4paper]{article}
\usepackage{latexsym}
\usepackage{graphicx}
\usepackage{amsmath}
\usepackage{amssymb}
\usepackage{epsfig}

\makeatletter \@addtoreset{equation}{section}
 \makeatother

\begin{document}

\title{ \textbf{ $\gamma^{*}$-semi-open Sets in Topological Spaces-II}}

\author{Bashir $Ahmad^{*}$ and Sabir $Hussain^{**}$ \\
{*} Centre for Advanced Studies in Pure and Applied Mathematics,
Bahauddin Zakariya \\ University, Multan, Pakistan.\\ {\bf Present
Address:} Department of Mathematics,\\ King Abdul Aziz University,
P. O. Box 80203, Jeddah 21589, Saudi Arabia.\\ E. mail:
drbashir9@gmail.com
\\{**} Department of Mathematics, Islamia University Bahawalpur,
Pakistan.\\ {\bf Present Address:} Department of Mathematics,\\
Yanbu University, P. O. Box 31387, Yanbu,  Saudi Arabia.\\ E.
mail: sabiriub@yahoo.com}
\date{}
\maketitle
%\maketitle \abstract{
\textbf{Abstract.} In this paper, we continue studying the
properties of   $\gamma^{*}$-semi-open sets in topological spaces
introduced by S. Hussain, B. Ahmad and T. Noiri[8]. We also
introduce and discuss the  $\gamma^{*}$-semi-continuous functions
which generalize semi-continuous functions defined by N. Levine
[10].
\\
\\
AMS(2000)Subject Classification. Primary 54A05, 54A10, 54D10.
%\begin{eqnarray*}
%\lefteqn {x=j+g+g_{h}}\\&& \subseteq k+i+g_g
%\\&& \subseteq gfF+DSSSSSSSSS+
%\end{eqnarray*}
\section {Introduction}
     A. Csaszar [6,7] defined generalized open sets in generalized
topological spaces. In 1975, Maheshwari and Prasad [11] introduced
concepts of semi-$T_1$-spaces and semi-$R_0$-spaces. In 1979, S.
Kasahara [9] defined an operation $\alpha$  on topological spaces.
In 1992 (1993), B. Ahmad and F.U. Rehman [1, 13] introduced the
notions of  $\gamma$-interior, $\gamma$-boundary and
$\gamma$-exterior points in topological spaces. They also studied
properties and characterizations of ($\gamma$ ,$\beta$
)-continuous mappings introduced by H. Ogata [12]. In 1999 (2005),
S. Hussain and B. Ahmad introduced the concept of $\gamma_
0$-compact, $\gamma^*$-regular, $\gamma$-normal spaces and
explored their many interesting properties. They initiated the
concept of $\gamma^*$-semi-open sets , $\gamma^*$-semi-closed
sets, $\gamma^*$-semi-closure, $\gamma^*$-semi-interior points in
topological spaces [8] and discussed many of  their properties. In
2006, they also introduced $\Lambda_{s}^{\gamma}$-set and
$\Lambda^{s^{\gamma}}$-set by using $\gamma^*$-semi-open sets.
Moreover, they showed that  the concepts of g.$\Lambda_{s}$-set ,
g.$\Lambda^{s}$-set, semi-$T_1$ space and semi-$R_0$ space can be
generalized by replacing semi- open sets with $\gamma^*$-semi-open
sets for an arbitrary monotone operator $\gamma \in \Gamma(X)$.

      In this paper, we continue studying the properties of
$\gamma^*$-semi-open sets in topological spaces introduced by S.
Hussain, B. Ahmad  and T. Noiri [8]. We also introduce and discuss
$\gamma$-semi-continuous functions which generalize
semi-continuous functions defined by N. Levine [10].

      Hereafter, we shall write space in place of topological space in the sequel.
We recall some definitions and results used in this paper to make
it self contained.\\
\section {Preliminaries}
 {\bf Definition 2.1 [9].}  Let (X,$\tau$) be a
space. An operation $\gamma$ : $\tau\rightarrow$ P(X) is a
function from $\tau$ to the power set of X such that V $\subseteq
V^\gamma$ , for each V $\in\tau$, where $V^\gamma$ denotes the
value of $\gamma$ at V. The operations defined by $\gamma$(G) = G,
$\gamma$(G) = cl(G) and
$\gamma$(G) = intcl(G) are examples of operation $\gamma$.\\
{\bf Definition 2.2 [9].}  Let A be a subset of a space X. A point
x $\in$ A is said to be $\gamma$-interior point of A iff there
exists an open nbd N of x such that $ N^\gamma\subseteq$ A. The
set of all such points is denoted by $int_\gamma$(A). Thus
\begin{center}
        $int_\gamma$ (A) = $\{ x \in A : x \in N \in \tau$  and  $N^\gamma\subseteq A \} \subseteq A$.
\end{center}

    Note that A is $\gamma$-open [12] iff A =$int_\gamma$(A). A set A is
called $\gamma$- closed [12] iff X-A is $\gamma$-open.\\
{\bf Definition 2.3 [9].}  A point x $\in$ X is called a
$\gamma$-closure point of A $\subseteq$ X, if $U^\gamma\cap A \neq
\phi$, for each open nbd U of x. The set of all $\gamma$-closure
points of A is called $\gamma$-closure of A and is denoted by
$cl_\gamma$(A). A subset A of X is called $\gamma$-closed, if
$cl_\gamma(A)\subseteq$A. Note that $cl_\gamma$(A) is contained in
every $\gamma$-closed superset
of A.\\
{\bf Definition 2.4 [12].} An operation $\gamma$ on $\tau$ is said
be regular, if for any open nbds U,V of x $\in$ X, there exists an
open nbd W of x such that $U^\gamma \cap V^\gamma\supseteq
W^\gamma$.\\
{\bf Definition 2.5 [12].}  An operation $\gamma$ on $\tau$ is
said to be open, if for every open nbd U of each $x \in X$,there
exists $\gamma$-open set B such that $x\in B$ and
$U^\gamma\subseteq
B$.\\

\section {Properties of $\gamma^{*}$-semi-open Sets and $\gamma^{*}$-semi-closed Sets}
{\bf Definition 3.1 [8].} A subset A of a space (X,$\tau$ ) is
said to be a $\gamma^{*}$-semi-open set, if  there exists a
$\gamma$ -open set O such that $O \subseteq A \subseteq
cl_{\gamma} (O)$. The set of all $\gamma^{*}$-semi-open sets is
denoted by
$SO_{\gamma_{*}}(X)$.\\

                Clearly every $\gamma$-open set is $\gamma^{*}$-semi-open set. The converse
may not be true as the following example shows:\\
{\bf Example 3.2.} Let X= $\{a,b,c\}$, $\tau =\{\phi , X, \{a\},
\{b\},\{a,b\},\{a,c\}\}$ be a topology on X. For $b\in X$, define
an operation $\gamma:\tau \rightarrow P(X)$ by
\begin{center}
$\gamma(A) = A^{\gamma}= \left\{
\begin{array}{ccc}
cl(A), & \mbox{if $b \in A$}\\
A, & \mbox{if $b \not\in A$}\\
\end{array} \right.$\\
\end{center}
Calculations gives $\tau_{\gamma}$ = $\{\phi , X, \{a\},
\{a,c\}\}$. Let A = $\{a,b\}$.Then there exists $\gamma$-open set
$\{a\}$ such that $\{a\}\subseteq A\subseteq cl_{\gamma} \{a\} =
X$. Also $SO_{\gamma}(X) = \{\phi , X, \{a\},
\{a,b\},\{a,c\}\}$.Thus A is
 $\gamma^{*}$-semi-open set but not $\gamma$-open set.\\

The following examples show that the concept of semi-open and
$\gamma^{*}$-semi-open sets are independent.\\
{\bf Example 3.3.} Let X= $\{a,b,c\}$, $\tau =\{\phi , X, \{a\},
\{b\},\{a,b\},\{a,c\}\}$ be a topology on X. For $b\in X$, define
an operation $\gamma:\tau \rightarrow P(X)$ by
\begin{center}
$\gamma(A) = A^{\gamma}= \left\{
\begin{array}{ccc}
A, & \mbox{if $a \in A$}\\
cl(A), & \mbox{if $a \not\in A$}\\
\end{array} \right.$\\
\end{center}
Calculations give $\gamma$-open sets are $\phi ,
X,\{a\},\{b\},\{a,b\}$. Then $A = \{b,c\}$ is
$\gamma^{*}$-semi-open but not a semi-open set.\\
{\bf Example 3.4.}  Let X= $\{a,b,c\}$, $\tau =\{\phi , X, \{a\},
\{b\},\{a,b\},\{a,c\}\}$ be a topology on X. For $b\in X$, define
an operation $\gamma:\tau \rightarrow P(X)$ by
\begin{center}
$\gamma(A) = A^{\gamma}= \left\{
\begin{array}{ccc}
A, & \mbox{if $b \in A$}\\
cl(A), & \mbox{if $b \not\in A$}\\
\end{array} \right.$\\
\end{center}
Calculations give $\gamma$-open sets are $\phi , X, \{b\},
\{a,b\},\{a,c\}$. Then $A = \{a\}$ is
semi-open but not  $\gamma^{*}$-semi-open.\\
{\bf Definition 3.5 [8].} A subset A of a space X is
$\gamma^{*}$-semi-closed
if and only if X-A is  $\gamma^{*}$-semi-open.\\
{\bf Definition 3.6 [8].} Let A be a subset of space X and $\gamma
\in \Gamma(X)$. The intersection of  all  $\gamma^{*}$-semi-closed
sets containing A is called $\gamma^{*}$-semi-closure of A and is
denoted by $scl_{\gamma^{*}}(A)$.\\
Note that A is
$\gamma^{*}$-semi-closed if and only if $scl_{\gamma^{*}}(A)=A $ .\\
{\bf Definition 3.7.} Let A be a subset of a space X and  $\gamma
\in \Gamma(X)$. The union of all $\gamma^{*}$-semi-open sets of X
contained in A is called
$\gamma^{*}$-semi-interior of A and is denoted by $sint_{\gamma^{*}}(A)$.\\
{\bf Definition 3.8.} A subset A of a space X is called
$\gamma^{*}$-semi-regular if and only if it is both
$\gamma^{*}$-semi-open and $\gamma^{*}$-semi-closed.

    If  A is  $\gamma^{*}$-semi-regular, then X-A is also  $\gamma^{*}$-semi-regular.
The class of all  $\gamma^{*}$-semi-regular sets of X is denoted
by $SR_{\gamma^{*}}(A)$.\\
{\bf Definition 3.9.}  A subset A of a space X is called
$\gamma^{*}$-pre-open, if  $A = int_{\gamma}(cl_{\gamma}(A))$.The
class of all $\gamma^{*}$-pre-open sets of X is denoted
by $PO_{\gamma^{*}}(A)$ .\\
{\bf Definition 3.10.}  Let A be a subset of a space X. Then the
$\gamma^{*}$-semi-boundary of A is denoted by
$sbd_{\gamma^{*}}(A)$ and is defined as $sbd_{\gamma^{*}}(A)=
scl_{\gamma^{*}}(A) \cap scl_{\gamma^{*}}(X-A)$.

    It is clear that $sbd_{\gamma^{*}}(A)$  is a $\gamma^{*}$-semi-closed set and contains $sbd(A)$, where
 $sbd(A)= scl(A) \cap scl(X-A)$\\
{\bf Remark 3.11.} It follows from the definition that
$sbd_{\gamma^{*}}(A)= \phi$  if and only if A is
$\gamma^{*}$-semi-regular.\\
{\bf Remark 3.12.}  Clearly from the definition of
$scl_{\gamma^{*}}(A)$  and $sint_{\gamma^{*}}(A)$ that if $A
\subseteq B$ then $scl_{\gamma^{*}}(A) \subseteq
scl_{\gamma^{*}}(B)$ and $sint_{\gamma^{*}}(A)\subseteq
sint_{\gamma^{*}}(B)$ .

      The following is immediate:\\
{\bf Theorem 3.13.}  Let A and B be any two subsets of a space X.
Then $scl_{\gamma^{*}}(A \cup B)=
 scl_{\gamma^{*}}(A) \cup scl_{\gamma^{*}}(B)$, where $\gamma$ is a regular
 operation.

          The following are straightforward:\\
{\bf Theorem 3.14.}  Let A be any subset of a space X. Then\\
(1)  $sint_{\gamma^{*}}(X-A) = X - scl_{\gamma^{*}}(A)$ .\\
(2)  $scl_{\gamma^{*}}(X-A) = X - sint_{\gamma^{*}}(A)$ .\\
(3)  $sint_{\gamma^{*}}(A) = X - scl_{\gamma^{*}}(X-A)$ .\\
{\bf Definition 3.15.}  The  $\gamma^{*}$-semi-exterior of A,
written $sext_{\gamma^{*}}(A)$, is
defined as the  $\gamma^{*}$-semi-interior of $X-A$. That is, $sext_{\gamma^{*}}(A) = sint_{\gamma^{*}}(X-A)$.\\
{\bf Theorem 3.16.}  In any space X, the following are
equivalent:\\
    (1)   X - $sbd_{\gamma^{*}}(A) = sint_{\gamma^{*}}(A)\cup
    sint_{\gamma^{*}}(X-A)$.\\
(2) $scl_{\gamma^{*}}(A) = sint_{\gamma^{*}}(A) \cup
sbd_{\gamma^{*}}(A)$.\\
(3) $sbd_{\gamma^{*}}(A) = scl_{\gamma^{*}}(A)\cap
scl_{\gamma^{*}}(X-A) = scl_{\gamma^{*}}(A)- sint_{\gamma^{*}}(A)$.\\
{\bf Proof.} $(1)\Rightarrow (2)$.  From (1), we have, $ X
-sbd_{\gamma^{*}}(A) = sint_{\gamma^{*}}(A)\cup
sint_{\gamma^{*}}(X-A)$. Also $X - sint_{\gamma^{*}}(X-A) =
scl_{\gamma^{*}}(A)$. Therefore $scl_{\gamma^{*}}(A) = X -
sint_{\gamma^{*}}(A)$.\\
$(2)\Rightarrow (3)$.  From (2), we have $scl_{\gamma^{*}}(A) =
sint_{\gamma^{*}}(A)\cup sbd_{\gamma^{*}}(A)=
scl_{\gamma^{*}}(A)\cap (X- sint_{\gamma^{*}}(A))$\\
$= scl_{\gamma^{*}}(A) \cap scl_{\gamma^{*}}(X-A)=
sbd_{\gamma^{*}}(A)$.\\
$(3)\Rightarrow (1)$. $sint_{\gamma^{*}}(A)\cup
sint_{\gamma^{*}}(X-A)=((sint_{\gamma^{*}}(A))^{c})^{c}\cup
(sint_{\gamma^{*}}(X-A))$\\
$=(scl_{\gamma^{*}}(X-A)\cap scl_{\gamma^{*}}(A))^{c}=
(sbd_{\gamma^{*}}(A))^{c} = X -
sbd_{\gamma^{*}}(A)$.\\
{\bf Proposition 3.17.}  (a)  From (3), we have
$sbd_{\gamma^{*}}(A) = sbd_{\gamma^{*}}(X-A)$.\\
 (b)$scl_{\gamma^{*}}(scl_{\gamma^{*}}(A))= scl_{\gamma^{*}}(A)$
, if  $\gamma$  is an open operation.\\
{\bf Theorem 3.18.}  Let A be any subset of a space X. Then\\
    (1)     $sbd_{\gamma^{*}}(A)=
    scl_{\gamma^{*}}(A)-sint_{\gamma^{*}}(A)$.\\
    (2)     $sbd_{\gamma^{*}}(A)\cap sint_{\gamma^{*}}(A) =\phi$
    .\\
    (3)     $scl_{\gamma^{*}}(A)=sint_{\gamma^{*}}(A)\cup
    sbd_{\gamma^{*}}(A)$.\\
    (4)     $sbd_{\gamma^{*}}(sint_{\gamma^{*}}(A))\subseteq
    sbd_{\gamma^{*}}(A)$.\\
    (5)   $sbd_{\gamma^{*}}(scl_{\gamma^{*}}(A))\subseteq sbd_{\gamma^{*}}(A)$ , if $\gamma$ is an open
    operation.\\
    (6)  $X-sbd_{\gamma^{*}}(A)= sint_{\gamma^{*}}(A)\cup sint_{\gamma^{*}}(X-A)$
    .\\
    (7)     $X = sint_{\gamma^{*}}(A)\cup sint_{\gamma^{*}}(X - A) \cup
    sbd_{\gamma^{*}}(A)$.\\
{\bf Proof.} (1)  $X- scl_{\gamma^{*}}(X-A) =
sint_{\gamma^{*}}(A)$ . Therefore $sbd_{\gamma^{*}}(A)=
scl_{\gamma^{*}}(A) \cap scl_{\gamma^{*}}(X-A)=
scl_{\gamma^{*}}(A)- (X - scl_{\gamma^{*}}(X-A))=
scl_{\gamma^{*}}(A) - sint_{\gamma^{*}}(A)$.\\
    (2) $sbd_{\gamma^{*}}(A)\cap sint_{\gamma^{*}}(A)=
(scl_{\gamma^{*}}(A)- sint_{\gamma^{*}}(A))\cap
sint_{\gamma^{*}}(A)= \phi$.\\
    (3) By (1), $sint_{\gamma^{*}}(A) \cup sbd_{\gamma^{*}}(A) =
    sint_{\gamma^{*}}(A)\cup (scl_{\gamma^{*}}(A) -
    sint_{\gamma^{*}}(A))= scl_{\gamma^{*}}(A)$.\\
    (4) and  (5)  follows directly.\\
    (6) $X - sbd_{\gamma^{*}}(A)= X - (scl_{\gamma^{*}}(A)\cap
    scl_{\gamma^{*}}(X-A))=(X - scl_{\gamma^{*}}(A))\cup (X -
    scl_{\gamma^{*}}(X-A))=sint_{\gamma^{*}}(X-A)\cup
    sint_{\gamma^{*}}(A) = sint_{\gamma^{*}}(A)\cup
    sint_{\gamma^{*}}(X-A)$.\\
(7)  is clear. Hence the proof.\\
{\bf Theorem 3.19.}  Let A be a subset of a space X. Then\\
            (1)     A is  $\gamma^{*}$-semi-open if and only if $A \cap sbd_{\gamma^{*}}(A)=\phi$
            .\\
            (2)     A is  $\gamma^{*}$-semi-closed if and only if $sbd_{\gamma^{*}}(A)\subseteq A$.\\
{\bf Proof.}  (1)  Let A be  $\gamma^{*}$-semi-open. Then  $X-A$
is $\gamma^{*}$-semi-closed and therefore
$scl_{\gamma^{*}}(X-A)=X-A$. Next (From [8]) $A \cap
sbd_{\gamma^{*}}(A)=A \cap (scl_{\gamma^{*}}(A)\cap
scl_{\gamma^{*}}(X-A))= A \cap scl_{\gamma^{*}}(A\cap
(X-A))=\phi$.

        Conversely, let $A\cap sbd_{\gamma^{*}}(A)=\phi$. Then $A\cap(scl_{\gamma^{*}}(A)\cap
scl_{\gamma^{*}}(X-A))= \phi$ , or $A\cap scl_{\gamma^{*}}(X-A)=
\phi$ implies $scl_{\gamma^{*}}(X-A) \subseteq (X-A)$. So $X-A$ is
$\gamma^{*}$-semi-closed implies A is  $\gamma^{*}$-semi-open.
This completes the proof.\\
(2)    Let A be  $\gamma^{*}$-semi-closed. Then
 $scl_{\gamma^{*}}(A)=A$. Now $sbd_{\gamma^{*}}(A)=
 scl_{\gamma^{*}}(A)\cap scl_{\gamma^{*}}(X-A)\subseteq
 scl_{\gamma^{*}}(A)=A$.

        Conversely, let $sbd_{\gamma^{*}}(A)\subseteq A$. Then $sbd_{\gamma^{*}}(A) \cap (X-A) = \phi$.
 Since $sbd_{\gamma^{*}}(A)= sbd_{\gamma^{*}}(X-A)$  implies that $sbd_{\gamma^{*}}(X-A) \cap (X-A) = \phi$.
Therefore $X-A$  is  $\gamma^{*}$-semi-open by (1). Hence A is
$\gamma^{*}$-semi-closed. This completes the proof.\\
{\bf Theorem 3.20.}  Let A be any  $\gamma^{*}$-semi-closed subset
of X. Then we have
$sbd_{\gamma^{*}}(sbd_{\gamma^{*}}(sbd_{\gamma^{*}}(A)))=
sbd_{\gamma^{*}}(sbd_{\gamma^{*}}(A))$, where $\gamma$ is a regular operation.\\
{\bf Proof.}  Since $scl_{\gamma^{*}}(A)=A$  if and only if A is
$\gamma^{*}$-semi-closed. We have\\
$sbd_{\gamma^{*}}(sbd_{\gamma^{*}}(sbd_{\gamma^{*}}(A)))=
scl_{\gamma^{*}}[sbd_{\gamma^{*}}(sbd_{\gamma^{*}}(A))] \cap
scl_{\gamma^{*}}[ X- sbd_{\gamma^{*}}(sbd_{\gamma^{*}}(A))]$\\
$= sbd_{\gamma^{*}}(sbd_{\gamma^{*}}(A))\cap scl_{\gamma^{*}}[ X
- (sbd_{\gamma^{*}}(sbd_{\gamma^{*}}(A)))]$  $ $ $ $ $ $ $ $ $ $ $ $ $.$ $.$ $.$   (*)\\
Now consider,\\
$X- (sbd_{\gamma^{*}}(sbd_{\gamma^{*}}(A)))=X - [
scl_{\gamma^{*}}(sbd_{\gamma^{*}}(A)) \cap scl_{\gamma^{*}}(X -
sbd_{\gamma^{*}}(A))] = X - [sbd_{\gamma^{*}}(A)\cap
scl_{\gamma^{*}}(X - sbd_{\gamma^{*}}(A))] = (X -
sbd_{\gamma^{*}}(A)) \cup (X - scl_{\gamma^{*}}(X -
sbd_{\gamma^{*}}(A)))$.\\
Therefore,\\
$scl_{\gamma^{*}}[ X - (sbd_{\gamma^{*}}(sbd_{\gamma^{*}}(A)))]=
scl_{\gamma^{*}}[ (X - sbd_{\gamma^{*}}(A)) \cup (X -
scl_{\gamma^{*}}(X - sbd_{\gamma^{*}}(A)))]$\\
$= scl_{\gamma^{*}}(X - sbd_{\gamma^{*}}(A)) \cup
scl_{\gamma^{*}}(X - scl_{\gamma^{*}}(X - sbd_{\gamma^{*}}(A)))$\\
$= D \cup scl_{\gamma^{*}}(X - scl_{\gamma^{*}}(X - D)) =
        X$ $ $ $ $ $ $ $ $ $ $ $.$ $.$ $.$  (**)\\
Where $D = scl_{\gamma^{*}}(X - sbd_{\gamma^{*}}(A))$. By (*) and
(**), we have\\
$sbd_{\gamma^{*}}(sbd_{\gamma^{*}}(sbd_{\gamma^{*}}(A)))=
sbd_{\gamma^{*}}(sbd_{\gamma^{*}}(A))\cap X =
sbd_{\gamma^{*}}(sbd_{\gamma^{*}}(A))$
 . Hence the proof.\\
{\bf Definition 3.21.} An operation $\gamma$  is said to be
semi-regular, if for every semi-open sets U and V containing $x
\in X$, there exists a
semi-open set W containing x such that $U^{\gamma} \cup V^{\gamma}\supseteq W^{\gamma}$ .\\
{\bf Definition 3.22.}  An operation  $\gamma$ is said to be
semi-open, if for every semi-open set U containing $x \in X$ ,
there exists a semi-open
set S containing x and $S \subseteq U^ { \gamma}$ .\\
{\bf Example 3.23.} Let $X = \{a,b,c\}$, $\tau = \{\phi  ,
X,\{a\},\{a,b\}\}$ be a topology on X. Define an
operation\\
$\gamma: \tau \rightarrow P(X)$ by $\gamma (A) =
cl(B)$. Clearly semi-open sets are $\{a\}, \{a,b\}, \{a,c\}, \phi,
X$ . Calculations give that $\gamma$ is not semi-regular but
$\gamma$ is semi-open.

    The following is straightforward:\\
{\bf Theorem 3.24.}  Let $\gamma: \tau \rightarrow P(X)$  be a
semi-regular operation on $\gamma$. If A and B are
 $\gamma^{*}$-semi-open, then $A \cap B$  is also  $\gamma^{*}$-semi-open.

          The following example shows that the condition that $\gamma$  is semi-regular is
necessary for the above theorem.\\
{\bf Example 3.25.}  $X = \{a,b,c\}$, $\tau = \{\phi  ,
X,\{a\},\{c\},\{a,c\}\}$ be a topology on X. Define an
operation\\
$\gamma: \tau \rightarrow P(X)$ by $\gamma (A) = int(cl(B))$.
$\gamma^{*}$-semi-open sets are  $\{a\}, \{c\}, \{a,b\}, \{b,c\},
\phi, X $ . It is easy to see that $\gamma$ is not semi-regular
operation. Here $\{a,b\}$ and $\{b,c\}$ are $\gamma^{*}$-semi-open
sets
but $\{a,b\} \cap \{b,c\} = \{b\}$ is not  $\gamma^{*}$-semi-open.\\
{\bf Theorem 3.26.}  Let A and B be subsets of a space X and
$\gamma$ a
semi-regular operation. Then\\
            (1)   $ sint_{\gamma^{*}}( sint_{\gamma^{*}}(A)) = sint_{\gamma^{*}}(A)$ .\\
            (2)   $ sint_{\gamma^{*}}(A \cup B) \supseteq sint_{\gamma^{*}}(A) \cup sint_{\gamma^{*}}(B) $.\\
            (3)    $ sint_{\gamma^{*}}(A \cap B)= sint_{\gamma^{*}}(A) \cap sint_{\gamma^{*}}(B) $.\\
{\bf Proof.}  (1) and (2) are obvious.\\
(3).  Let $x \in sint_{\gamma^{*}}(A) \cap sint_{\gamma^{*}}(B)$.
Then $ x \in sint_{\gamma^{*}}(A)$ and $x \in
sint_{\gamma^{*}}(B)$ . This gives that there exist semi-open sets
U and V containing x such that $U^{\gamma} \subseteq A$  and $ V ^
{\gamma} \subseteq B$ . Therefore, we have $U^{\gamma} \cap
V^{\gamma} \subseteq A \cap B$. Since $\gamma$ is semi-regular,
there exists semi-open nbd W of x such that$ W^{\gamma} \subseteq
U^{\gamma} \cap V^{\gamma}$. Consequently, $W^{\gamma} \subseteq A
\cap B$. That is, $ x \in sint_{\gamma^{*}}(A \cap B)$ implies
$sint_{\gamma^{*}}(A) \cap sint_{\gamma^{*}}(B) \subseteq
sint_{\gamma^{*}}(A \cap B)$.

        Conversely, $ A \cap B \subseteq A$ implies $sint_{\gamma^{*}}(A \cap B) \subseteq sint_{\gamma^{*}}(A)$ .
Also $ A \cap B \subseteq B$ , gives $sint_{\gamma^{*}}(A \cap B)
\subseteq sint_{\gamma^{*}}(B)$ . Thus $ sint_{\gamma^{*}}(A \cap
B) \subseteq sint_{\gamma^{*}}(A) \cap sint_{\gamma^{*}}(B) $ .
Consequently, $ sint_{\gamma^{*}}(A \cap B)= sint_{\gamma^{*}}(A) \cap sint_{\gamma^{*}}(B) $ .
This completes the proof.\\
{\bf Theorem 3.27.}  For any subsets A, B of  a space X , if   is
semi-regular operation, then\\
            (1) $ sext_{\gamma^{*}}(A \cup B) = sext_{\gamma^{*}}(A) \cap sext_{\gamma^{*}}(B) $    .\\
            (2)  $sbd_{\gamma^{*}}(A \cup B) = [ sbd_{\gamma^{*}}(A) \cap scl_{\gamma^{*}}(X - B)] \cup [ sbd_{\gamma^{*}}(B)
            \cap scl_{\gamma^{*}}(X - A)]$.\\
            (3) $sbd_{\gamma^{*}}(A \cap B) = [ sbd_{\gamma^{*}}(A) \cap scl_{\gamma^{*}}(B)] \cup [ sbd_{\gamma^{*}}(B)
            \cap scl_{\gamma^{*}}(A)]$.\\
{\bf Proof.}  (1) $ sext_{\gamma^{*}}(A \cup B)=
sint_{\gamma^{*}}(X - (A \cup B)) = sint_{\gamma^{*}}((X - A) \cap
(X - B))$\\
$= sint_{\gamma^{*}}(X - A) \cap sint_{\gamma^{*}}(X - B)=
sext_{\gamma^{*}}(A)\cap sext_{\gamma^{*}}(B)$.\\
  (2) $sbd_{\gamma^{*}}(A \cup B) = scl_{\gamma^{*}}( A \cup B)
  \cap scl_{\gamma^{*}}(X - (A \cup B))$\\
$ = (scl_{\gamma^{*}}(A) \cup scl_{\gamma^{*}}(B))\cap
scl_{\gamma^{*}}((X -A) \cap (X - B))$\\
$= (scl_{\gamma^{*}}(A) \cup scl_{\gamma^{*}} (B))\cap
(scl_{\gamma^{*}}(X-A)\cap scl_{\gamma^{*}}(X - B))$\\
$=(scl_{\gamma^{*}}(A) \cup scl_{\gamma^{*}}(B)) \cap
(scl_{\gamma^{*}}(X -A)\cap scl_{\gamma^{*}}(X -B))$\\
$=(scl_{\gamma^{*}}(A)\cap scl_{\gamma^{*}}(X -A))\cap
(scl_{\gamma^{*}}(X -B) \cup scl_{\gamma^{*}}(B)) \cap
(scl_{\gamma^{*}}(X-A)\cap scl_{\gamma^{*}}(X -B))$\\
$=sbd_{\gamma^{*}}(A) \cap scl_{\gamma^{*}}(X -B)\cup
sbd_{\gamma^{*}}(B)\cap scl_{\gamma^{*}}(X -A).$\\
 (3) $ sbd_{\gamma^{*}}(A \cap B) = scl_{\gamma^{*}}(A \cap B)
 \cap scl_{\gamma^{*}}(X - (A \cap B))$\\
$= (scl_{\gamma^{*}}(A) \cap scl_{\gamma^{*}}(B)) \cap
scl_{\gamma^{*}}((X -A) \cup (X -B))$\\
$= (scl_{\gamma^{*}}(A) \cup scl_{\gamma^{*}}(B)) \cap
(scl_{\gamma^{*}}((X -A) \cap scl_{\gamma^{*}}(X -B)))$\\
$= (scl_{\gamma^{*}}(A) \cap scl_{\gamma^{*}}(X -A)) \cap
(scl_{\gamma^{*}}(X -B) \cup scl_{\gamma^{*}}(B)) \cap
(scl_{\gamma^{*}}(X -A)\cap scl_{\gamma^{*}}(X -B))$\\
$= sbd_{\gamma^{*}}(A) \cap scl_{\gamma^{*}}(X -B) \cup
sbd_{\gamma^{*}}(B) \cap scl_{\gamma^{*}}(X -A)$.\\
{\bf Proposition 3.28.}  If  $\gamma$  is semi regular operation, then\\
                (1)  $ sext_{\gamma^{*}}( X - sext_{\gamma^{*}} (A)) = sext_{\gamma^{*}} (A)$.\\
                (2)  $ sext_{\gamma^{*}}(A \cap B) \supseteq sext_{\gamma^{*}}(A) \cup sext_{\gamma^{*}}(B)$.

    The following is immediate from the definition:\\
{\bf Proposition 3.29.}  A set A of a space X is said to be
$\gamma^{*}$-semi-closed if and only if there exists a
$\gamma$-closed set F such
that  $int_{\gamma}(F) \subseteq A \subseteq F$.\\
{\bf Theorem 3.30.}  A subset A of a space X is
$\gamma^{*}$-semi-closed if and only if $int_{\gamma}(cl_{\gamma}(A)) \subseteq A$.\\
{\bf Proof.}  Suppose that A is  $\gamma^{*}$-semi-closed. Then
there exists a $\gamma$-closed set F such that $int_{\gamma}(F)
\subseteq A \subseteq F$. This implies that
$int_{\gamma}(cl_{\gamma}(A)) \subseteq int_{\gamma} (F) \subseteq
A$ gives $int_{\gamma}(cl_{\gamma} (A)) \subseteq A$.\\
Next suppose that $int_{\gamma}(cl_{\gamma} (A)) \subseteq A$ .
Putting $cl_{\gamma}(A)= F$, we get
 $int_{\gamma}(F)
\subseteq A \subseteq cl_{\gamma}(A) = F$. This shows that A is
$\gamma$-semi-closed. Hence the proof.\\

\section {$\gamma$-semi-continuous  and $\gamma$-semi-open functions }
{\bf Definition 4.1.}  A function $f:(X, \tau)\rightarrow (Y,
\tau$) is a $\gamma$-semi-continuous if and
only if for any  $\gamma$-open B in Y, $f^{-1}(B)$  is  $\gamma^{*}$-semi-open in X.\\
{\bf Theorem 4.2.} Let $f:X \rightarrow Y $  be a function and $x
\in X$ . Then $f$  is $\gamma$-semi-continuous if and only if for
each $\gamma$-open set B containing f(x), there exists $A \in
SO_{\gamma^{*}}(X)$ such that $x \in A$ and $f(A) \subseteq B$,
where $\gamma$ is a regular
operation.\\
{\bf  Proof.} Suppose that $f$ is  $\gamma$-semi-continuous.
Therefore for $\gamma$-open set B in Y, $f^{-1}(B)$  is
 $\gamma^{*}$-semi-open in X. We prove that for each  $\gamma$-open set B containing $f(x)$, there exists
 $\gamma^{*}$-semi-open set A in X such that $x \in A$  and $f(A) \subseteq B$ . Let $x \in f^{-1}(B) \in SO_{\gamma^{*}}(X)$ and
 $A=f^{-1}(B)$  . Then $x \in A$  and $f(A) \subseteq ff^{-1}(B) \subseteq B$, where  B is  $\gamma$-open.

         Conversely, let B be  $\gamma$-open set in Y. We prove that inverse image of  $\gamma$-open set in Y is
$\gamma^{*}$-semi-open set in X.  Let $x \in f^{-1}(B)$ . Then
$f(x) \in B$. Thus there exists an  $A_{x} \in SO_{\gamma^{*}}(X)$
such that  $x \in A_{x}$ and $f(A_{x}) \subseteq B$ . Then $x \in
A_{x} \subseteq f^{-1}(B)$ and  $f^{-1}(B) = \bigcup_{x \in f^{-1}
(B)} A_{x}$ implies $f^{-1}(B) \in SO_{\gamma^{*}}(X)$, since
$\gamma$ is regular operation [8]. This proves that f is
$\gamma$-semi-continuous.\\
{\bf Definition 4.3.} A function  $f:X \rightarrow Y $ is said to
be  $\gamma$-semi-open if and only if for each
 $\gamma$-open set U in X, $f(U)$  is  $\gamma^{*}$-semi-open in Y.\\
{\bf Lemma 4.4.}  If A is  $\gamma^{*}$-semi-open and $A \subseteq
B$. Then $A \subseteq cl_{\gamma}(int_{\gamma}(B))$.\\
{\bf Proof.} Since A is $\gamma^{*}$-semi-open, $A \subseteq
cl_{\gamma}(int_{\gamma}(A))$[8]. Also, $A \subseteq B$ implies
that $int_{\gamma}(A) \subseteq int_{\gamma} (B)$ . So
$cl_{\gamma}(int_{\gamma}(A)) \subseteq cl_{\gamma}( int_{\gamma}
(B))$ implies $A \subseteq cl_{\gamma}( int_{\gamma}
(B))$. This completes the proof.\\
{\bf Theorem 4.5.}  A function $f:X \rightarrow Y $ is
$\gamma$-semi-open if and only if for every subset\\
$E \subseteq X$ , we have $f(int_{\gamma}(E)) \subseteq
cl_{\gamma}(int_{\gamma}(f(E))).$\\
{\bf  Proof.}  Let f be $\gamma^{*}$-semi-open. Since
$f(int_{\gamma}(E)) \subseteq f(E)$, then by supposition
$f(int_{\gamma}(E))$ is $\gamma$-semi-open. So by Lemma 4.4, we
get $f(int_{\gamma}(E)) \subseteq
cl_{\gamma}(int_{\gamma}(f(E))).$

     Conversely, let the condition hold and let G be any  $\gamma$-open set in X. Then\\
 $int_{\gamma}(f(G))\subseteq f(G)
\subseteq f( int_{\gamma}(G)) \subseteq
cl_{\gamma}(int_{\gamma}(f(G)))$. Therefore $f(G)$ is
$\gamma^{*}$-semi-open and consequently $f$ is $\gamma$-semi-open.
This proves the theorem.

       Similarly, we can prove:\\
{\bf Theorem 4.6.}  A function $f:X \rightarrow Y $ is
$\gamma$-semi-open if and only if for every
subset G  of Y, $int_{\gamma}(f^{-1}(G)) \subseteq cl_{\gamma}( f^{-1}(int_{\gamma}(G))) $.\\
{\bf Theorem 4.7.}  Let  $f:X \rightarrow Y $  be
$\gamma$-semi-continuous and  $\gamma$-semi-open and let $A \in
SO_{\gamma^{*}}(X)$. Then $f(A) \in SO_{\gamma^{*}}(Y)$, where
$\gamma$ is an open operation.\\
{\bf Proof.}  Since A is $\gamma$-semi-open, there exists
$\gamma$-open set O in X such that $O \subseteq A \subseteq
cl_{\gamma} (O)$. Therefore $f(O) \subseteq f(A) \subseteq
f(cl_{\gamma} (O)) \subseteq cl_{\gamma}(f(O))$.
Thus $f(A) \in SO_{\gamma^{*}}(X)$. Hence the proof.\\
{\bf Theorem 4.8.}  Let $f:X \rightarrow Y $   be a function. Then
the following are
equivalent:\\
            (1)    f is  $\gamma$-semi-continuous.\\
            (2)    $f(scl_{\gamma^{*}}(A)) \subseteq cl_{\gamma}(f(A))$ , for any subset A of
            X.\\
            (3)     $sbd_{\gamma^{*}}(f^{-1}(B))\subseteq f^{-1}(bd_{\gamma}(B))$, for any subset B of Y.\\
{\bf Proof.}$(1)\Rightarrow (2)$. Let A be any subset of X. Since
$cl_{\gamma}(f(A))$ is $\gamma$-closed in Y, then f  is
$\gamma$-semi-continuous implies $f^{-1}(cl_{\gamma}(f(A)))$ is
$\gamma^{*}$- semi-closed in X, which contains A. $A \subseteq
f^{-1}(cl_{\gamma}(f(A))$ gives\\
$scl_{\gamma^{*}}(A) \subseteq
scl_{\gamma^{*}}(f^{-1}(cl_{\gamma}(f(A)))=
f^{-1}(cl_{\gamma}(f(A))$.
 Therefore  $f(scl_{\gamma^{*}}(A)) \subseteq
 f(f^{-1}(cl_{\gamma}(f(A))))$. Consequently,
$f(scl_{\gamma^{*}}(A)) \subseteq cl_{\gamma}(f(A))$ . This gives
(2).\\
$(2)\Rightarrow (1)$. Suppose $f(scl_{\gamma^{*}}(A)) \subseteq
cl_{\gamma}(f(A))$ , for $A \subseteq X$ . To prove (1), let B be
$\gamma$-closed subset of Y.
 We show that  $f^{-1}(B)$ is  $\gamma^{*}$-semi-closed. By our hypothesis,\\
 $f(scl_{\gamma^{*}}(f^{-1}(B))) \subseteq
cl_{\gamma}(f(f^{-1}(B))) \subseteq cl_{\gamma}(B) = B$ implies\\
$scl_{\gamma^{*}}(f^{-1}(B)) \subseteq
f^{-1}f(scl_{\gamma^{*}}(f^{-1}(B))) \subseteq f^{-1}(B)$ or
$scl_{\gamma^{*}}(f^{-1}(B)) \subseteq
 f^{-1}(B)$
 implies $f^{-1}(B)$ is  $\gamma^{*}$-semi-closed in X. Thus f is
 $\gamma$-semi-continuous.\\
 $(1)\Rightarrow (3)$.  Suppose that f is $\gamma$-semi-continuous. Let B be any subset of Y.
 Then\\
 $sbd_{\gamma^{*}}(f^{-1}(B))= scl_{\gamma^{*}}(f^{-1}(B))\cap scl_{\gamma^{*}}(X - f^{-1}(B))$
 \\
$ \subseteq scl_{\gamma^{*}}(f^{-1}(cl_{\gamma}(B))\cap
scl_{\gamma^{*}}(f^{-1}(cl_{\gamma}(Y - B))$\\
$ = f^{-1} (cl_{\gamma}(B)) \cap f^{-1} (cl_{\gamma}(Y - B)) =
f^{-1} (cl_{\gamma}(B)\cap cl_{\gamma} (Y - B)) =
f^{-1}(bd_{\gamma}(B)).$\\
Therefore, $sbd_{\gamma^{*}}(f^{-1}(B))\subseteq
f^{-1}(bd_{\gamma}(B))$.\\
$(3)\Rightarrow (1)$ . Let B be  $\gamma$-closed in Y. We show
that $ f^{-1}(B) \in scl_{\gamma^{*}}(X)$. Since
$sbd_{\gamma^{*}}(f^{-1}(B)) \subseteq
f^{-1}(bd_{\gamma}(B)\subseteq f^{-1}(B))$ implies
$sbd_{\gamma^{*}}(f^{-1}(B))\subseteq f^{-1}(B))$  gives
$f^{-1}(B)$ is $\gamma^{*}$-semi-closed in X. Hence $f$ is
$\gamma$-semi-continuous.\\
{\bf Theorem 4.9.}  A function $f:X \rightarrow Y $ is
$\gamma$-semi-continuous if and only
if for every subset G  of Y, $ scl_{\gamma^{*}}(f^{-1}(G)) \subseteq cl_{\gamma}(f^{-1}(G))$ .\\
{\bf Proof.}  Since  $f$ is $\gamma$-semi-continuous, so by
Theorem 4.8, we have
\begin {center}
$f(scl_{\gamma^{*}}(A)) \subseteq cl_{\gamma}(f(A))$ $ $ $ $ $ $ $
$ $.$ $.$ $.$  (1)
\end {center}
Let G  be any subset of  Y. Put $A = f^{-1}(G)$  then
$f(scl_{\gamma^{*}}(f^{-1}(G))) \subseteq cl_{\gamma}(ff^{-1}(A))
\subseteq cl_{\gamma}(G)$ . This implies that
$f(scl_{\gamma^{*}}(f^{-1}(G)))\subseteq cl_{\gamma}(G)$.

       Conversely, suppose that for for any subset G  of Y, $f(scl_{\gamma^{*}}(f^{-1}(G)))\subseteq cl_{\gamma}(G)$ .
 Let  $G = f(A)$
for $A \subseteq X$. Then $scl_{\gamma^{*}}(A) \subseteq
scl_{\gamma^{*}}(f^{-1}(G)) \subseteq f^{-1}(cl_{\gamma}(f(A))$.
Therefore $f(scl_{\gamma^{*}}(A)) \subseteq cl_{\gamma}(f(A))$. By
Theorem 4.8, $f$
is  $\gamma$-semi-continuous.\\
{\bf Lemma 4.10.}  Let $A \subseteq B$  and $B \in
SC_{\gamma^{*}}(X)$ , then $sbd_{\gamma^{*}}(A) \subseteq B$,
where $SC_{\gamma^{*}}(X)$ is the class of
all $\gamma^{*}$-semi-closed sets of X.\\
{\bf Proposition 4.11.}  If $A \cap B = \phi$ and A is
$\gamma$-open, then $A \cap cl_{\gamma}(B) = \phi$.

           Using this fact, we obtain the following.\\
{\bf Lemma 4.12.}  Let $A \cap B = \phi$  and A a  $\gamma$-open set of X, then $A \cap bd_{\gamma}(B) = \phi$ .\\
{\bf Theorem 4.13.}  Let B be any subset of a space Y. A bijective
function $f:X \rightarrow Y $ is
 $\gamma$-semi-open if and only if  $f^{-1}(sbd_{\gamma^{*}}(B) \subseteq bd_{\gamma}(f^{-1}(B))$ .\\
{\bf Proof.} Suppose f is $\gamma$-semi-open and  $B \subseteq Y$
. Put
\begin {center}
$ U = X - bd_{\gamma}(f^{-1}(B))$$ $ $ $ $ $ $ $ $ $ $.$ $.$ $.$
(*)
\end {center}
Then U is $\gamma$-open. Therefore $f(U)$  is
$\gamma^{*}$-semi-open in Y. This gives $Y - f(U) \in
SC_{\gamma^{*}}(Y)$ . Since f is bijective, therefore by (*), we
get $B \subseteq Y - f(U)$. Then by Lemma 4.10, we have\\
$f^{-1}(sd_{\gamma^{*}}(B)) \subseteq f^{-1}(Y - f(U)) = f^{-1}(Y)
- f^{-1}(f(U)) \subseteq X -U$\\
$= X - (X - bd_{\gamma}(f^{-1}(B))) = bd_{\gamma}(f^{-1}(B))$  Consequently, we have\\
$f^{-1}(sd_{\gamma^{*}}(B)) \subseteq bd_{\gamma}(f^{-1}(B))$ .

              Conversely, suppose $f^{-1}(sd_{\gamma^{*}}(B)) \subseteq bd_{\gamma}(f^{-1}(B))$ , for any subset B of Y.
We prove that f is $\gamma$-semi-open. Let U be $\gamma$-open in
X. Put $B = Y - f(U)$. Clearly $B \cap f(U) = \phi$, which gives
 $U \cap f^{-1}(B) = \phi$. By Lemma 4.12, we have $U \cap bd_{\gamma}(f^{-1}(B)) = \phi$. Consequently,
$f^{-1}(sd_{\gamma^{*}}(B)) \subseteq bd_{\gamma}(f^{-1}(B))$
gives $U \cap f^{-1}(sbd_{\gamma^{*}}(B)) = \phi$ or $\phi = f( U
\cap f^{-1} ( sbd_{\gamma^{*}}(B))) = f(U) \cap
sbd_{\gamma^{*}}(B)$ implies $ sbd_{\gamma^{*}}(B) \subseteq Y -
f(U) = B$ gives B is $\gamma^{*}$-semi-closed. Consequently $f(U)$
is $\gamma^{*}$-semi-open. This completes the
proof.\\
{\bf Theorem 4.14.} A function  $f:X \rightarrow Y $ is
$\gamma$-semi-open if and only if for every subset $A \subseteq X$
, $f(int_{\gamma}(A)) \subseteq cl_{\gamma}(int_{\gamma}f(A))$.\\
Proof. Since $int_{\gamma}(A) \subseteq A$ and let $f$ be
$\gamma$-semi-open, we have $f(int_{\gamma}(A)) \subseteq f(A)$.
By supposition $f(int_{\gamma}(A))$ is $\gamma^{*}$-semi-open. So
by Lemma 4.4, we get $ f(int_{\gamma}(A)) \subseteq
cl_{\gamma}(int_{\gamma}(f(A)))$.

           Conversely, let G be any $\gamma$-open set in X. Then $f(G) = f( int_{\gamma}(G))\subseteq
cl_{\gamma}(int_{\gamma}(f(G)))$ .Therefore $f(G)$
 is  $\gamma^{*}$-semi-open and consequently f is $\gamma$-semi-open. This completes the
 proof.\\

\end{document}